\documentclass [a4paper,twoside,13pts]{article}
\usepackage [francais]{babel}
\usepackage [latin1] {inputenc}
\usepackage [T1] {fontenc}
\usepackage{amsmath}
\usepackage{stmaryrd}
\usepackage{graphicx}
\title{Pénalisations de l'araignée brownienne (Penalizations of Walsh's
Brownian motion) }
\author {Joseph Najnudel}
\begin{document} 
\maketitle 
\noindent
\textbf{Résumé : }
 Dans cet article, nous pénalisons une araignée brownienne $(A_t)_{t
 \geq 0}$ sur un ensemble fini $E$ de demi-droites, avec un poids égal
 à $Z \exp (\alpha_{N_t}
 X_t + \gamma L_t)$, où $t$ est un nombre positif, $(\alpha_k)_{k \in
 E}$ une famille de réels indexée par $E$, $\gamma$ un paramètre réel, $X_
t$ la distance de $A_t$ à l'origine, $N_t \in E$ la demi-droite associée à
 $A_t$, $L_t$ le temps local de $(A_s)_{0 \leq s \leq t}$ à l'origine, et $Z$ une
constante de normalisation. Nous montrons que la famille des mesures
 de probabilité obtenue par ces pénalisations converge vers une mesure
 limite quand $t$ tend vers l'infini, et nous étudions le comportement
 de cette mesure limite. \\ \\
\textbf{Abstract : } In this paper, we penalize a Walsh's Brownian
 motion $(A_t)_{t \geq 0}$ on a finite set $E$ of rays, with a weight
 equal to $Z \exp (\alpha_{N_t} X_t + \gamma L_t)$, where $t$ is a
 positive number, $(\alpha_k)_{k \in E}$ a family of real numbers
 indexed by $E$, $\gamma$ a real parameter, $X_t$ the distance from
 $A_t$ to the origin, $N_t \in E$ the ray associated to $A_t$, $L_t$ the
 local time of $(A_s)_{0 \leq s \leq t}$ at the origin, and $Z$ a
 constant of normalization. We show that the family of the probability measures
 obtained by these penalizations converges to a limit measure as $t$
 tends to infinity, and
 we study the behaviour of this limit measure. \\ \\
\textbf{Mots-clé : } pénalisation, temps local, araignée
 brownienne. 
\textbf{classifications AMS : } 60B10, 60J65 (60G17, 60G44, 60J25,
 60J55).  
\section{Présentation du problème et des résultats principaux obtenus}
\noindent
Récemment, de nombreuses études de pénalisations du
mouvement brownien ont été effectuées, en particulier par B. Roynette, P. Vallois et
M. Yor (voir [6], [7], [8]). \\
Dans [8], les pénalisations étudiées sont fonctions de la valeur $X_t$
atteinte par un mouvement brownien en un temps $t$, et du suprémum
$S_t$ de ce mouvement
brownien jusqu'au temps $t$. Plus précisément, on considère une
famille de mesures de probablité $(\mathbf{W}^{(t)})_{t \geq 0}$
sur $\mathcal{C} (\mathbf{R}_+, \mathbf{R})$ vérifiant, pour tout
$\Gamma_t$ appartenant à la tribu $\mathcal{F}_t$ engendrée
par $(X_s)_{s \in [0,t]}$ ($(X_t)_{t \geq 0}$ étant le processus
canonique de $\mathcal{C} (\mathbf{R}_+, \mathbf{R})$) :
$$\mathbf{W}^{(t)}(\Gamma_t) = \frac{\mathbf{W} [\mathbf{1}_{\Gamma_t}
  f(X_t,S_t)]}{\mathbf{W} [f(X_t,S_t)]} $$ où 
$S_t$ est le maximum de $X_s$ pour $s \in [0,t]$, $\mathbf{W}$ la
mesure de Wiener, et $f$ une fonction de $\mathbf{R}^2$ dans
$\mathbf{R}_+$. \\ 
B. Roynette, P. Vallois et M. Yor montrent alors que pour certains
  choix de la fonction $f$, il existe une mesure de probabilité
  $\mathbf{W}^{(\infty)}$ sur $\mathcal{C} (\mathbf{R}_+, \mathbf{R})$
  telle que pour tout $s \geq 0$ et tout $\Gamma_s \in \mathcal{F}_s$
  : $$ \mathbf{W}^{(t)} (\Gamma_s) \underset{t \rightarrow
  \infty}{\rightarrow} \mathbf{W}^{(\infty)} (\Gamma_s)$$
Un des cas où cette convergence a lieu est celui où $f(a,y)= \exp(\lambda y + \mu a)$
  avec $\lambda$, $\mu \in \mathbf{R}$. \\ \\  
Par un changement de mouvement brownien, les résultats de [8]
peuvent alors être adaptés au cas où $S_t$ est remplacé par $L_t$ (temps
local en 0 de $(X_u)_{u \leq t}$), et où $X_t$ est remplacé par $L_t -
|X_t|$; en effet, d'après le théorème d'équivalence de Lévy : $(S_t -
X_t, S_t)_{t \geq 0}$ a même loi que $(|X_t|,L_t)_{t \geq 0}$. \\ \\
Dans ces conditions, les poids exponentiels étudiés dans [8] prennent
la forme : $Z \exp (\alpha |X_t| + \gamma L_t)$ où $\alpha$ et $\gamma$
sont des paramètres réels et où $Z$ est une constante de
normalisation. \\ \\
Le but de notre article est de généraliser l'étude de ces
pénalisations exponentielles à toutes les araignées browniennes dont
les trajectoires se situent sur un nombre fini de demi-droites
concourantes (voir [1] et [9] pour une description de ces processus).
\\ \\ \\ 
Soit $(E,\mu)$ un espace de probabilité fini; on suppose $\mu(\{ m \} ) >
0$ pour tout $m \in E$. \\
On considère, sur l'espace $\mathbf{R}_E = \{ (0,0) \} \cup
(\mathbf{R}_+^* \times E)$, la distance $d$ définie par :
$$d((x,k),(y,l)) = |x-y| \mathbf{1}_{k = l} + (x+y) \mathbf{1}_{k \neq l}$$  
Cette distance permet de considérer $\mathcal{C}_E$, espace des
fonctions continues de $\mathbf{R}_+$ dans  $\mathbf{R}_E$, et de
munir cet espace de la tribu $\mathcal{T}_E$ de la convergence
uniforme. \\ \\
$(A_t = (X_t,N_t))_{t \geq 0}$ désigne alors le processus canonique (à
valeurs dans $\mathbf{R}_E$) associé à l'espace $(\mathcal{C}_E,
\mathcal{T}_E)$ et on note, pour tout $t \in \mathbf{R}_+$,
$\mathcal{F}_t$ la sous-tribu de $\mathcal{T}_E$ engendrée par
$(A_s)_{0 \leq s \leq t}$, $A_s$ étant donc une application de
$(\mathcal{C}_E, \mathcal{T}_E)$ dans $(\mathbf{R}_E,d)$ (muni de sa
tribu borélienne). \\ \\
Pour $(x,k) \in \mathbf{R}_E$, on peut alors considérer, sur
$\mathcal{C}_E$, la mesure de probabilité $\mathbf{W}_{(E,\mu,x,k)}$, sous
laquelle $(A_t)_{t \geq 0}$ est une araignée brownienne à trajectoires
dans $\mathbf{R}_E$, issue de $(x,k)$, et telle que pour tous $s$, $t$,
$s \leq t$, la loi de $N_t$ sachant $X_s = 0$ est $\mu$. \\
Rappelons (voir [1]) que cette araignée brownienne est un processus de
Feller qu'il est possible de caractériser par son semi-groupe
$(P_t)_{t \geq 0}$; pour toute fonction $f$ borélienne
bornée : $$ P_t f(x,k) = 2 \underset{m \in E}{\sum} \mu_m
\underset{\mathbf{R}_+^*}{\int} dy p_t(x+y) f(y,m) +
\underset{\mathbf{R}_+^*}{\int} dy (p_t(x-y) - p_t(x+y)) f(y,k)$$ 
avec $\mu_m = \mu( \{ m \} )$ et $p_t(a) = \frac{1}{\sqrt{2 \pi t}}
e^{-a^2/2t}$. \\ \\
Sous $\mathbf{W}_{(E, \mu, 0,0)}$, le processus $(X_t)_{t \geq 0}$ est un
mouvement  brownien réfléchi (comme on le voit directement avec le
semi-groupe). \\
Soit $\mathcal{I}$ l'ensemble des intervalles d'excursion de $(X_t)_{t
  \geq 0}$, $N_t$ est alors constant sur chaque intervalle $I \in
\mathcal{I}$ : on peut donc poser $N_t = N_I$ pour $t \in I$. On
montre alors que conditionnellement à $(X_t)_{t \geq 0}$, les
$(N_I)_{I \in \mathcal{I}}$ sont des variables aléatoires
indépendantes et de loi $\mu$. \\ \\
A présent, nous allons effectuer un changement de probabilité par
rapport à l'araignée brownienne, de la manière suivante : pour $\alpha
= (\alpha_i)_{i \in E}$ une famille de réels indexée par $E$, $\gamma
\in \mathbf{R}$ et $t \in \mathbf{R}_{+}$, on pose  
$$ \mathbf{W}_{(E, \mu)}^{(t, \alpha, \gamma)} = \frac{\exp(\alpha_{N_t} X_t +
  \gamma L_t)}{\mathbf{W}_{(E, \mu, 0,0)} [\exp(\alpha_{N_t} X_t +
  \gamma L_t)]} . \mathbf{W}_{(E, \mu, 0,0)}$$
où $L_t$ est le temps local de $(X_t)_{t \geq 0}$ : 
$$ L_t = \underset{\epsilon \rightarrow 0}{\lim \inf} \frac{1}{2
  \epsilon} \int_0^t \mathbf{1}_{X_s \leq \epsilon} ds$$ (en fait, la
  limite inférieure ci-dessus est presque sûrement une limite). 
Dans ces conditions, nous allons prouver les deux théorèmes suivants : \\
  \\
\textbf{Théorème 1 : } Il existe une mesure de probabilité $\mathbf{W}
_{(E,\mu)}^{(\infty, \alpha,
  \gamma)}$ telle que pour tout $s \in \mathbf{R}_+$ et tout $\Gamma_s
  \in \mathcal{F}_s$ : $$
  \mathbf{W}_{(E,\mu)}^{(t,\alpha,\gamma)}(\Gamma_s) \underset{t \rightarrow
  \infty}{\rightarrow} \mathbf{W}_{(E,\mu)}^{(\infty, \alpha, \gamma)}
  (\Gamma_s)$$ 
De plus, pour tout $s \in \mathbf{R}_+$, la restriction de
  $\mathbf{W}_{(E,\mu)}^{(\infty, \alpha, \gamma)}$ à $\mathcal{F}_s$
  est équivalente à la loi de l'araignée brownienne
  sur $[0,s]$ associée à $(E, \mu)$ (nous en préciserons la densité
  dans la suite de cet article). \\ \\
\textbf{Remarque : } On voit clairement que si le théorème 1 est vrai,
  la famille des densités de $\mathbf{W}_{(E,\mu)}^{(\infty, \alpha,
  \gamma)}$ par rapport à $\mathbf{W}_{(E,\mu,0,0)}$,
  conditionnellement à $(\mathcal{F}_s)_{s \geq 0}$, est une
  $(\mathcal{F}_s)_{s \geq 0}$-martingale sous la probabilité
  $\mathbf{W}_{(E,\mu,0,0)}$. \\ \\ 
\textbf{Théorème 2 : } Le processus canonique $(A_t)_{t \geq
  0}$ sous la mesure $\mathbf{W}_{(E,\mu)}^{(\infty, \alpha, \gamma)}$ peut
  être décrit de la manière suivante : \\ \\
- Si $\gamma \geq \alpha_m$ pour tout $m$ et $\gamma > 0$,
  $(X_s)_{s \geq 0}$ est la valeur absolue d'un processus bang-bang de
  paramètre $\gamma$, ou encore : $X_t = S_t - Y_t$, où
  $(Y_t)_{t \geq 0}$ est un mouvement brownien avec drift $\gamma$ et
  $S_t$ son suprémum sur $[0,t]$; et
  $(N_s)_{s \geq 0}$ est construit à partir de $(X_s)_{s \geq 0}$ de
  la même manière que pour l'araignée brownienne : chaque excursion est
  sur la branche $m$ avec probabilité $\mu_m$, indépendamment des
  autres excursions (voir également [3]). \\ \\
- Si $\bar{\alpha} = \max (\alpha) > \gamma$
  et $\bar{\alpha} > 0$, $(X_s)_{s \geq 0}$ est un processus dont
  la loi a une densité égale à $\frac{\bar{\alpha} -
    \gamma}{\bar{\alpha}} \exp(\gamma L_{\infty})$ par rapport à la loi
  de la valeur absolue d'un mouvement brownien avec drift
  $\bar{\alpha}$ (dont $L_{\infty}$ est le temps local sur
  $\mathbf{R}_+$), et $(N_s)_{s \geq 0}$ est obtenu en effectuant la même
  démarche que pour l'araignée initiale, puis en conditionnant le
  résultat par le fait que la dernière excursion de $(X_s, N_s)_{s
  \geq 0}$ se situe sur une branche $m$ vérifiant $\alpha_m =
  \bar{\alpha}$. \\ \\
- Si $\gamma = 0$ et $\alpha_m \leq 0$ pour tout $m$, 
  $(X_s,N_s)_{s \geq 0}$ est une araignée brownienne associée à $(E, \mu)$. \\ \\
- Si $\gamma < 0$ et $\alpha_m \leq 0$ pour tout $m$, on considère $(Y_s,R_s)_{s \geq 0}$ une araignée brownienne,
  $\mathbf{e}$ une variable exponentielle de paramètre $\gamma$ indépendante de
  $(Y_s,R_s)_{s \geq 0}$, $\tau_{\mathbf{e}}$ l'inverse du temps local
  de $(Y_s)_{s \geq 0}$ en $\mathbf{e}$, $(Z_s)_{s \geq 0}$ un
  processus de Bessel de dimension 3 issu de 0 et indépendant des
  variables précédentes, et $M$ une variable aléatoire (également
  indépendante des précédentes) définie sur $E$. \\
Si la loi de $M$ est bien choisie (voir la section 4 pour plus de
  détails sur cette loi), le processus $(X_s,N_s)_{s \geq 0}$ a même loi que
  $(\tilde{X}_s, \tilde{N}_s)_{s \geq 0}$, avec $(\tilde{X}_s, \tilde{N}_s) = (Y_s,R_s)$ pour $s \leq
  \tau_{\mathbf{e}}$, et  $(\tilde{X}_s, \tilde{N}_s) = (Z_{s-\tau_{\mathbf{e}}}, M)$ pour $s
  \geq \tau_{\mathbf{e}}$. \\ \\ \\
Dans la suite de cet article, nous allons tout d'abord évaluer
  l'expression : $\mathbf{W}_{(E, \mu, x,k)} [\exp(\alpha_{N_t} X_t + \gamma
  L_t)]$, puis nous utiliserons cette évaluation pour démontrer les deux
  théorèmes annoncés.
 
\section{Etude de l'expression $\mathbf{W}_{(E, \mu, x,k)} [\exp(\alpha_{N_t}
  X_t + \gamma L_t)]$}
Afin de prouver l'existence de $\mathbf{W}_{(E,\mu)}^{(\infty, \alpha,
  \gamma)}$, nous allons commencer par chercher une expression
qui majore $\mathbf{W}_{(E, \mu, x,k)} [\exp(\alpha_{N_t} X_t + \gamma
  L_t)]$ tout en étant équivalente à cette quantité quand $t$ tend
vers l'infini. \\ \\
Pour cela, nous allons tout d'abord effectuer une étude de deux
quantités, $I(\beta,
\gamma, x,t)$ et $J( \beta, x,t)$ ($\beta, \gamma \in \mathbf{R}$, $x,
t \in \mathbf{R}_+$), définies de la manière suivante :
 $$ I(\beta, \gamma, x,t) = \mathbf{E}_x [\exp(\beta |Y_t| + \gamma
  L_t) \mathbf{1}_{T \leq t}]$$
$$J(\beta, x,t) =
\mathbf{E}_x [\exp(\beta Y_t) \mathbf{1}_{T > t}]$$
 où $(Y_t)_{t \geq
  0}$ est un mouvement brownien issu de $x$, $L_t$ est le temps
local en zéro de $(Y_s)_{0 \leq s \leq t}$, et $T = \inf \{s \geq 0,
Y_s = 0 \}$. \\ \\
Comme on le voit dans les expressions définissant $I$ et $J$, l'étude
de ces quantités ne fait intervenir que des propriétés du mouvement
brownien. \\ \\ 
 \textbf{Etude de $J(\beta,x,t)$ : } Le principe de réflexion
 implique : $$ J(\beta,x,t) = \mathbf{E}_x [e^{\beta Y_t}
 \mathbf{1}_{T > t}] = 
\mathbf{E}_x [e^{\beta Y_t}
 \mathbf{1}_{Y_t > 0}] - 
 \mathbf{E}_x [e^{\beta Y_t}
 \mathbf{1}_{Y_t>0, T \leq t}]$$ $$ = \mathbf{E}_x
 [e^{\beta Y_t} \mathbf{1}_{Y_t > 0}] - 
\mathbf{E}_x [e^{- \beta Y_t} \mathbf{1}_{Y_t < 0}] $$ $$=
\frac{1}{\sqrt{2 \pi t}} \int_0^\infty (e^{- ((x-y)^2/2t) +
    \beta y} - 
e^{-((x+y)^2/2t) + \beta y}) dy$$
\textbf{Supposons $\beta <0$ : } On a la majoration immédiate : 
$$ e^{-(x-y)^2/2t} - e^{-(x+y)^2/2t} \leq \frac{(x+y)^2}{2t} -
\frac{(x-y)^2}{2t} = \frac{2xy}{t}$$
On en déduit : 
$$ J(\beta,x,t) \leq \sqrt{\frac{2}{\pi t^3}} x \int_0^\infty y
e^{\beta y} dy = \sqrt{\frac{2}{\pi t^3}} \frac{x}{\beta^2}$$
Par ailleurs, on a les encadrements suivants : 
$$ 1 - \frac{(x-y)^2}{2t} \leq e^{-(x-y)^2/2t} \leq 1 -
\frac{(x-y)^2}{2t} + \frac{(x-y)^4}{8t^2}$$ 
$$ 1 - \frac{(x+y)^2}{2t} \leq e^{-(x+y)^2/2t} \leq 1 -
\frac{(x+y)^2}{2t} + \frac{(x+y)^4}{8t^2}$$
La différence entre $J(\beta,x,t)$ et son majorant donné précédemment
est donc bornée par : $$ \frac{1}{\sqrt{2 \pi t}} \int_0^\infty
\frac{(x+y)^4}{8 t^2} e^{\beta y} dy = C(x) t^{-5/2}$$
où $C(x)$ ne dépend que de $x$. \\ \\ 
On en déduit que : $$L(\beta,x,t) = \sqrt{\frac{2}{\pi t ^3}}
\frac{x}{\beta^2}$$ est à la fois un majorant et un équivalent de
$J(\beta,x,t)$ quand $t \rightarrow \infty$ ($x$ étant fixé). \\ \\
\textbf{Supposons $\beta = 0$ : } On obtient ici $$ J(\beta,x,t) =
\frac{1}{\sqrt{2 \pi t}} \int_{-x}^x e^{-y^2/2t} dy$$ expression
admettant comme majorant et comme équivalent : $$ L(\beta,x,t) =
\sqrt{\frac{2}{\pi t}} x$$ 
\textbf{Supposons $\beta > 0$ : } On a l'égalité suivante : 
$$ \frac{1}{\sqrt{2 \pi t}} \int_{-\infty}^{\infty} (e^{
  -((x-y)^2/2t) + \beta y} -e^{-((x+y)^2/2t) + \beta y}) dy $$ $$ =  \frac{1}{\sqrt{2 \pi
  t}} \int_{-\infty}^{\infty} dz e^{-(z^2/2t) + \beta z} (e^{\beta x} -
  e^{- \beta x}) = 2 \sinh (\beta x) e^{t \beta ^2/2}$$ 
Or : $$  \frac{1}{\sqrt{2 \pi t}} \int_{-\infty}^{0} (e^{-((x-y)^2/2t) +
  \beta y}-e^{-((x+y)^2/2t) + \beta y}) dy = -J(- \beta,x,t)$$
D'où l'égalité : 
$$ J(\beta,x,t) = J(-\beta,x,t) + 2 \sinh(\beta x) e^{t \beta^2/2}$$
On a donc le majorant et équivalent suivant : $$ L(\beta,x,t) = 2
\sinh(\beta x) e^{t \beta^2/2} + \sqrt{\frac{2}{\pi t^3}}
\frac{x}{\beta^2}$$ 
\\ 
On peut alors regrouper tous les cas possibles grâce à l'expression
suivante : $$ L(\beta,x,t) = \sqrt{\frac{2}{\pi t^3}}
\frac{x}{\beta^2} \mathbf{1}_{\beta \neq 0} + \sqrt{\frac{2}{\pi t}} x
\mathbf{1}_{\beta=0} + 2 \sinh(\beta x) \exp(t \beta^2/2)
\mathbf{1}_{\beta > 0}$$ 
\\
\textbf{Etude de la loi de $(|Y_t|,L_t)$} \\ \\
Avant de procéder à l'évaluation de $I(\beta, \gamma, x,t)$, nous
allons étudier la loi jointe de $(|Y_t|,L_t)$, lorsque $(Y_t)_{t \geq
  0}$ est un mouvement brownien issu de $x$ et $(L_t)_{t \geq 0}$ son
temps local en zéro. \\ 
Plus précisément, nous allons prouver le lemme suivant : \\ \\
\textbf{Lemme : } Avec les hypothèses précédentes, on a les deux
résultats suivants : \\ \\
- Pour $z \in \mathbf{R}_{+}$, $\mathbf{P} [L_t+|Y_t| \in dz,
  L_t > 0] = \sqrt{\frac{2}{\pi t^3}} z(x+z) \exp \left( -
    \frac{(x+z)^2}{2t} \right) dz$. \\ \\
- Conditionnellement au fait que $L_t> 0$, $\Theta_t =
\frac{|Y_t|}{L_t + |Y_t|}$ est une variable uniforme sur $[0,1]$,
indépendante de $L_t + |Y_t|$. \\ \\
\textbf{Preuve : } Dans le calcul suivant, $\mathbf{P}_y$ désigne la
loi d'un mouvement brownien issu de $y$, $(Y_t)_
{t \geq 0}$ est le processus canonique de $\mathcal{C}
(\mathbf{R}_+, \mathbf{R})$, $L_t$ son temps local, et $T
= \inf \{ t \geq 0, Y_t = 0 \}$. \\ 
On a, pour $y \in \mathbf{R}_+$ et $l>0$ :  
$$ \mathbf{P}_x (|Y_t| \in dy, L_t \in dl) = \underset{s_1+ s_2 \leq
  t}{\int} \mathbf{P}_x (T \in ds_1) \mathbf{P}_0(|Y_{t-s_1}| \in dy)
  $$ $$\mathbf{P}_0 (\underset{0 \leq u \leq t-s_1}{\sup} \{u | Y_u =
 0 \} \in d(t-s_2), L_{t-s_1} \in dl | |Y_{t-s_1}|=y)$$ 
Par renversement du temps effectué sur le pont brownien :  
$$ \mathbf{P}_x (|Y_t| \in dy, L_t \in dl) $$ $$= \underset{s_1+s_2
  \leq t}{\int} \mathbf{P}_x (T \in ds_1)  \mathbf{P}_y (T \in ds_2,
  L_{t-s_1} \in dl, |Y_{t-s_1}| \in [0,dy])$$ $$ =  \underset{s_1+s_2
  \leq t}{\int} \mathbf{P}_x (T \in ds_1) \mathbf{P}_y (T \in ds_2)
  \mathbf{P}_0 (|Y_{t-s_1-s_2}| \in [0,dy], L_{t-s_1-s_2} \in dl)$$
  $$= \underset{s_1+s_2 \leq t}{\int} \frac{2 dy} {\sqrt{2 \pi
  (t-s_1-s_2)}}  \mathbf{P}_x (T \in ds_1) \mathbf{P}_y (T \in ds_2) 
 \mathbf{P}_0 (L_{t-s_1-s_2} \in dl | Y_{t-s_1-s_2} = 0)$$
Or la loi du temps local d'un pont brownien sur $[0,t-s_1-s_2]$ est
  connue : c'est la loi de la racine carrée d'une variable
  exponentielle dont le paramètre est $\frac{1}{2(t-s_1-s_2)}$. \\ 
On en déduit : 
$$ \mathbf{P}_x (|Y_t| \in dy, L_t \in dl) =  \underset{s_1 + s_2
  \leq t}{\int} dy dl \frac{2 l e^{-l^2/2(t-s_1-s_2)}}{\sqrt{2 \pi
  (t-s_1-s_2)^3}} \mathbf{P}_x (T \in ds_1) \mathbf{P}_y (T \in ds_2)
  $$ $$ = 2 dy dl \underset{s_1 + s_2 \leq t}{\int} \mathbf{P}_x(T \in
  ds_1) \mathbf{P}_y(T \in ds_2) \frac{\mathbf{P}_l (T \in
  d(t-s_1-s_2))}{d(t-s_1-s_2)} $$ $$ = 2 dy dl \underset{s_1 + s_2
  \leq t}{\int} ds_1 ds_2 D_x(s_1) D_y(s_2) D_l(t-s_1-s_2) $$ $$ = 2
  dy dl D_x * D_y * D_l (t) = 2 D_{x+y+l} (t) dy dl$$ $$ =
  \sqrt{\frac{2}{\pi t^3}} (l+x+y) \exp \left( - \frac{(l+x+y)^2}{2t} \right) dy dl$$
$D_a(u)$ désignant la densité de $T$ en $u$ sous $\mathbf{P}_a$. \\ \\
Ces égalités impliquent le lemme annoncé. \\ \\
Remarquons que ce lemme est fortement lié au théorème de Pitman (voir
  également [5]). \\ \\
\textbf{Etude de $I(\beta, \gamma,x,t)$} \\ \\
Le lemme précédent permet d'écrire la formule suivante : $$ I(
  \beta, \gamma, x,t) = \mathbf{E} \left[ \sqrt{\frac{2}{\pi t ^3}}
  \int_0^{\infty} z(x+z) \exp \left( \frac{-(x+z)^2}{2t}+ \Phi z \right) dz \right]$$ où
  $\Phi$ est une variable uniforme sur $[\beta, \gamma]$. \\ 
Distinguons à présent plusieurs cas. \\ \\
\textbf{Supposons $\beta$, $\gamma < 0$ : } Dans ce cas, le théorème
  de convergence monotone prouve que $\mathbf{E} [\int_0^{\infty}
  z(x+z) e^{-((x+z)^2/2t) + \Phi z} dz]$ croît vers $\mathbf{E}
  [\int_0^{\infty} z(x+z) e^{\Phi z} dz ]$ quand $t$ tend vers
  l'infini. \\ \\ 
Or pour $\phi \in \mathbf{R}_-^*$, $\int_0^{\infty} z(x+z) e^{\phi z}
  dz = \frac{x} {\phi^2} + \frac{2}{|\phi|^3}$.  \\ \\
On en déduit que si $\beta \neq \gamma$ : $$ \mathbf{E} \left[
  \int_0^{\infty} z (x+z) e^{\Phi z} dz \right] = \frac{1}{\gamma -
  \beta} \int_{\beta}^{\gamma} \left(\frac{x}{\phi^2} +
  \frac{2}{|\phi|^3} \right) d \phi $$ $$ = \frac{1}{\gamma - \beta}
  \left[-\frac{x}{\phi} + \frac{1}{\phi^2} \right]_{\beta}^{\gamma} =
  \frac{x}{\beta \gamma} + \frac{|\beta| +|\gamma|}{\beta^2
  \gamma^2}$$ et que cette dernière égalité se prolonge en fait au cas
  où $\beta = \gamma$. \\ \\
On en déduit que $I(\beta, \gamma, x,t)$ admet comme majorant et comme
  équivalent : $$ K( \beta, \gamma, x,t) = \sqrt{\frac{2}{\pi t^3}}
  \left(\frac{x}{\beta \gamma} + \frac{|\beta| + |\gamma|}{\beta^2
  \gamma^2} \right)$$ 
\textbf{Supposons $\beta = 0$, $\gamma < 0$ : } On a, pour tout $z$;
  $$\mathbf{E} [e^{\Phi z}] = \frac{1}{|\gamma|} \int_{\gamma}^0
  e^{\phi z} d \phi = \frac{1 - e^{\gamma z}}{|\gamma| z}$$ ($\Phi$
  étant une variable uniforme sur $[\gamma, 0]$). \\ D'où : $$I(\beta,
  \gamma, x,t) = \frac{1}{|\gamma|} \sqrt{\frac{2}{\pi t^3}}
  \int_0^{\infty} (x+z) e^{-(x+z)^2/2t} dz $$ $$- \frac{1}{|\gamma|}
  \sqrt{\frac{2}{\pi t^3}} \int_0^\infty (x+z) e^{- ((x+z)^2/2t) +
  \gamma z} dz$$
On a $\int_0^{\infty} (x+z) e^{\gamma z} < \infty$, donc le deuxième
  terme de l'expression ci-dessus, négatif, est dominé par $t^{-3/2}$
  quand $t$ tend vers l'infini. \\ \\
Par ailleurs, $$ \int_0^{\infty} (x+z) e^{-(x+z)^2/2t} dz = \left[ -t
  e^{-(x+z)^2/2t} \right]_0^{\infty} = t e^{-x^2/2t}$$ admet $t$ comme
  majorant et comme équivalent quand $t$ tend vers l'infini. \\ \\
Ces deux propriétés permettent d'en déduire que $I(\beta, \gamma,
  x,t)$ admet comme majorant et équivalent : $$ K(\beta, \gamma, x,t)
  = \frac{1}{|\gamma|} \sqrt{\frac{2}{\pi t}}$$
\textbf{Supposons $\gamma = 0$, $\beta < 0$ : } On a évidemment par
  symétrie : $$K(\beta, \gamma, x,t) = \frac{1}{|\beta|}
  \sqrt{\frac{2}{\pi t}}$$
\textbf{Supposons $\beta=\gamma= 0$ : } On a : $$ I(\beta, \gamma,
  x,t) = \sqrt{\frac{2}{\pi t^3}} \int_0^{\infty} z(x+z)
  e^{-(x+z)^2/2t} dz $$ $$ = \sqrt{\frac{2}{\pi t^3}} \left[-tz
  e^{-(x+z)^2/2t} \right]_0^{\infty} + \sqrt{\frac{2}{\pi t}}
  \int_0^{\infty} e^{-(x+z)^2/2t} dz = 2 \mathbf{P}(\mathcal{N} \geq
  x/\sqrt{t}) $$ où $\mathcal{N}$ est une variable gaussienne centrée
  réduite. \\ 
On peut donc prendre : $$ K(\beta, \gamma, x,t) = 1$$
\textbf{Supposons $\beta > 0$ et $\beta > \gamma$ : } On a : $$
  \mathbf{E}[e^{\Phi z} ] = \frac{1}{\beta - \gamma}
  \int_{\gamma}^{\beta} e^{\phi z} d \phi = \frac{e^{\beta z} -
  e^{\gamma z}}{(\beta - \gamma) z}$$ ($\Phi$ uniforme sur $[\gamma,
  \beta]$). \\ 
On en déduit : $$ I(\beta, \gamma, x,t)= \frac{1}{\beta - \gamma}
  \sqrt{\frac{2}{\pi t^3}} \int_0^{\infty} (x+z) e^{-((x+z)^2/2t) +
  \beta z} dz $$ $$- \frac{1}{\beta - \gamma} \sqrt{\frac{2}{\pi t^3}}
  \int_0^{\infty} (x+z) e^{-((x+z)^2/2t) + \gamma z} dz$$
Pour tout $\phi > 0$ : $$ \sqrt{\frac{2}{\pi t^3}} \int_0^{\infty}
  (x+z) e^{-((x+z)^2/2t) + \phi z} dz $$ $$= \sqrt{\frac{2}{\pi t^3}} \left[
  -t e^{-((x+z)^2/2t) + \phi z} \right]_0^{\infty} + \phi
  \sqrt{\frac{2}{\pi t} } \int_0^{\infty} e^{-((x+z)^2/2t) + \phi z} dz
  $$ $$ = \sqrt{\frac{2}{\pi t}} e^{-x^2/2t} + 2 \phi e^{- \phi x + t
  \phi ^2/2} - \phi \sqrt{\frac{2}{\pi t}} \int_0^{\infty}
  e^{-((x-z)^2/2t) - \phi z} dz$$ avec $$ \int_0^{\infty} e^{-((x-z)^2/2t)
  - \phi z} dz \leq \int_0^{\infty} e^{-\phi z} dz = \frac{1}{\phi}$$
Donc la quantité ci-dessus admet $ \sqrt{2/\pi t} + 2 \phi e^{-\phi x+
  t \phi^2/2}$ comme majorant et comme équivalent. \\
On peut en particulier en déduire que le second terme de $I(\beta,
  \gamma, x,t)$, négatif, est négligeable devant le premier quand $t$
  tend vers l'infini (quel que soit le signe de $\gamma$). \\ \\ 
Tout ceci permet de prendre : $$ K(\beta, \gamma, x, t) =
  \frac{1}{\beta- \gamma} \sqrt{\frac{2}{\pi t}} + \frac{2
  \beta}{\beta - \gamma} \exp(-\beta x + t \beta ^2/2)$$ 
\textbf{Supposons $\gamma > 0$ et $\gamma > \beta$ : } La symétrie
  permet d'obtenir immédiatement : $$ K(\beta,\gamma,x,t) =
  \frac{1}{\gamma-\beta} \sqrt{\frac{2}{\pi t}} + \frac{2
  \gamma}{\gamma-\beta} \exp(- \gamma x + t \gamma^2/2)$$
\textbf{Supposons $\gamma = \beta > 0$ : } On a ici : $$ I(\beta,
  \gamma,x,t) = \sqrt{\frac{2}{\pi t^3}} \int_0^{\infty } z(x+z)
  e^{-((x+z)^2/2t)  + \gamma z} dz$$ 
Or $$\int_0^{\infty} (z(x+z) - \gamma t z -t) e^{-((x+z)^2/2t) + \gamma
  z} dz = \left[ -t z e^{-((x+z)^2/2t) + \gamma z} \right]_0^{\infty} =
  0$$
On en déduit : $$ I(\beta,\gamma,x,t) = \sqrt{\frac{2}{\pi t}}
  \int_0^{\infty} (\gamma z + 1) e^{-((x+z)^2/2t) + \gamma z} dz$$ 
Par ailleurs, on a : $$ \sqrt{\frac{2}{\pi t}} \int_0^{\infty}
  e^{-((x+z)^2/2t) + \gamma z} dz = 2 e^{- \gamma x + t \gamma^2/2} -
  \sqrt{\frac{2}{\pi t}} \int_0^{\infty} e^{-((x-z)^2/2t) - \gamma z} dz$$
quantité équivalente et inférieure à $2 e^{-\gamma x + t
  \gamma^2/2}$.  \\ \\
La quantité $- \gamma x \sqrt{\frac{2}{\pi t}} \int_0^{\infty}
  e^{-((x+z)^2/2t) + \gamma z} dz$ est donc négative et équivalente à
  $-2 \gamma x e^{- \gamma x + t \gamma^2/2}$. \\ \\
De plus, d'après un calcul précédemment effectué, $\gamma
  \sqrt{\frac{2}{\pi t}} \int_0^{\infty} (x+z) e^{-((x+z)^2/2t) + \gamma
  z} dz$ est équivalent et inférieur à $\gamma \sqrt{\frac{2t}{\pi }}
  + 2 t \gamma^2 e^{-\gamma x + t \gamma ^2/2}$ (voir l'étude du cas
  $\beta > 0$ et $\beta > \gamma$). \\ \\
En additionnant les trois termes évalués ci-dessus, on obtient :
  $$K(\beta, \gamma, x,t) = \gamma \sqrt{\frac{2t}{\pi}} + 2(t
  \gamma^2 + 1) \exp(-\gamma x + t \gamma^2/2)$$
On a donc le tableau suivant : \\ \\ \\
\begin{tabular}{|p{3.3 cm}|p{4.5 cm}|p{3.2 cm}|} 
\hline
Conditions sur $\beta$ et $\gamma$ & $K(\beta, \gamma,x,t)$ &
Equivalent le plus simple quand $t \rightarrow \infty$ \\
\hline
$\beta$, $\gamma < 0$ & $\sqrt{\frac{2}{\pi t^3}} \left(\frac{x}{\beta
    \gamma} + \frac{|\beta|+|\gamma|}{\beta^2 \gamma ^2} \right)$ &
$\sqrt{\frac{2}{\pi t^3}} \left(\frac{x}{\beta \gamma} +
  \frac{|\beta|+|\gamma|}{\beta^2 \gamma^2} \right)$ \\
\hline
$\beta = 0$, $\gamma <0$ & $\frac{1}{|\gamma|} \sqrt{\frac{2}{\pi t}}$
& $\frac{1}{|\gamma|} \sqrt{\frac{2}{\pi t}}$ \\
\hline
$\gamma = 0$, $\beta < 0$ & $\frac{1}{|\beta|} \sqrt{\frac{2}{\pi t}}$
& $\frac{1}{|\beta|} \sqrt{\frac{2}{\pi t}}$ \\
\hline 
$\beta = \gamma = 0$ & $1$ & $1$ \\
\hline
$\beta > 0$, $\beta > \gamma$ & $\frac{1}{\beta - \gamma}
\sqrt{\frac{2}{\pi t}} + \frac{2 \beta}{\beta - \gamma} e^{-\beta x +
  t \beta^2/2}$ & $\frac{2 \beta}{\beta - \gamma} e^{- \beta x + t
  \beta^2/2}$ \\
\hline
$\gamma > 0$, $\gamma > \beta$ & $\frac{1}{\gamma - \beta}
\sqrt{\frac{2}{\pi t}} + \frac{2 \gamma}{\gamma - \beta} e^{-\gamma x
  + t \gamma^2/2}$ & $\frac{2 \gamma}{\gamma - \beta} e^{-\gamma x + t
  \gamma^2/2}$ \\
\hline 
$\gamma = \beta > 0$ & $\gamma \sqrt{\frac{2 t}{\pi}} + 2(t \gamma^2 +
1) e^{- \gamma x + t \gamma^2/2}$ & $2 t \gamma^2 e^{- \gamma x + t
  \gamma^2/2}$ \\
\hline 
 
\end{tabular}
\\ \\ \\
A présent, nous avons obtenu des majorants et des équivalents pour les
quantités $I$ et $J$ et nous sommes en mesure d'évaluer l'expression
$\mathbf{W}_{(E, \mu,x,k)} [\exp(\alpha_{N_t} X_t + \gamma L_t)]$. \\ \\
En effet, on a : $$\mathbf{W}_ {(E, \mu, x,k)}[\exp(\alpha_{N_t} X_t + \gamma
  L_t)] = A_1 + A_2$$ avec 
$$ A_1 = \mathbf{W}_{(E, \mu, x,k)} [\exp(\alpha_{N_t} X_t + \gamma L_t)
\mathbf{1}_{T \leq t}]$$
$$ A_2 = \mathbf{W}_{(E, \mu, x,k)} [\exp(\alpha_{N_t} X_t + \gamma L_t)
\mathbf{1}_{T > t}]$$
où $T = \inf \{ s \geq 0, X_s = 0 \}$. 
D'après la propriété de Markov de l'araignée, conditionnellement au
fait que $(X_s)_{s \geq 0}$ s'annule avant $t$, $N_t$ est une
variable de loi $\mu$, indépendante de $(X_t, L_t)$. D'où : 
$$A_1 = \underset{m \in E}{\sum} \mu_m I(\alpha_m, \gamma, x,t)$$
Par ailleurs, si $(X_s)_{s \geq 0}$ ne s'annule pas avant $t$, il est
évident que $L_t = 0$ et $N_t = k$. \\
On a donc $A_2 = J(\alpha_k,x,t)$, et on en déduit : $$
\mathbf{W}_ {(E, \mu, x,k)}[e^{\alpha_{N_t} X_t + \gamma
  L_t}] = \underset{m \in E} {\sum} \mu_m I(\alpha_m, \gamma,x,t) +
J(\alpha_k,x,t)$$
Cette égalité permet de prouver la proposition suivante : \\ \\
\textbf{Proposition : } La quantité $$ Q_{(E,\mu)}(\alpha, \gamma, x,k,t) = \underset{m \in
  E}{\sum} \mu_m K(\alpha_m, \gamma,x,t) + L(\alpha_k, x,t)$$ 
où $K(\alpha_m, \gamma, x,t)$ et $L( \alpha_k, x, t)$ ont été définis
précédemment dans cette section, vérifie les propriétés suivantes : $$
\mathbf{W}_ {(E, \mu, x,k)} [\exp(\alpha_{N_t} X_t + \gamma L_t)] \leq Q_{(E,\mu)}
(\alpha, \gamma, x,k,t)$$
$$ \mathbf{W}_{(E, \mu, x,k)} [\exp(\alpha_{N_t} X_t + \gamma L_t)]
\underset{t \rightarrow \infty}{\sim} Q_{(E,\mu)}
(\alpha, \gamma, x,k,t)$$
\section{Preuve de l'existence de la mesure $\mathbf{W}_{(E,\mu)}^{(\infty,
  \alpha, \gamma)}$}
On observe, tout d'abord, que quels que soient $\beta$ et $\gamma$, il
existe $C(\beta, \gamma)$ tel que pour tous $x$, $t$ : $$K(\beta,
\gamma, x,t) \leq C(\beta, \gamma) (1+x) K(\beta, \gamma,0,t)$$ (en
fait, $K(\beta, \gamma, x, t) \leq K(\beta, \gamma, 0,t)$ dès que
$\sup(\beta, \gamma) \geq 0$). \\ \\
On en déduit l'existence de $C_{(E,\mu)}(\alpha, \gamma)$, tel que pour tous
$t$, $x$ : $$\underset{m \in E}{\sum} \mu_m K(\alpha_m, \gamma,x,t)
\leq C_{(E,\mu)} (\alpha, \gamma) (1 + x) \underset{m \in E}{\sum} \mu_m
K(\alpha_m, \gamma, 0, t) $$ $$= C_{(E,\mu)}(\alpha, \gamma) (1+x)
Q_{(E,\mu)} (\alpha,
\gamma, 0,0,t)$$
A présent, fixons $\beta$ et $\gamma$ dans $\mathbf{R}$, et supposons
$t \geq 1$. \\
Si $\beta < 0$ et $\gamma < 0$, $L(\beta,x,t) = \sqrt{\frac{2}{\pi t
    ^3}} \frac{x}{\beta^2}$ et $ K(\beta, \gamma, 0,t) =
\sqrt{\frac{2}{\pi t^3}} \frac{|\beta|+|\gamma|}{\beta^2 \gamma^2}$,
ce qui implique : $$ L(\beta,x,t) = \frac{x \gamma^2}{|\beta| +
  |\gamma|} K(\beta, \gamma, 0,t)$$
Si $\beta<0$ et $\gamma=0$, $ K(\beta,\gamma,0,t) = \frac{1}{|\beta|}
\sqrt{\frac{2}{\pi t}} \geq \frac{1}{|\beta|} \sqrt{\frac{2}{\pi t
^3}} $, et donc : $$ L(\beta,x,t) \leq \frac{x}{|\beta|}
K(\beta,\gamma,0,t)$$
Si $\beta <0$ et $\gamma >0$, $K(\beta,\gamma,0,t) \geq
\frac{1}{\gamma - \beta} \sqrt{\frac{2}{\pi t}} \geq \frac{1}{\gamma -
  \beta} \sqrt{\frac{2}{\pi t^3}} $, d'où : $$ L(\beta,x,t) \leq x
\left(\frac{\gamma-\beta}{\beta^2} \right) K(\beta,\gamma,0,t)$$
Si $\beta = 0$ et $\gamma < 0$, $L(\beta,x,t) = \sqrt{\frac{2}{\pi t}}
x$ et $K(\beta, \gamma,0,t) = \frac{1}{|\gamma|} \sqrt{\frac{2}{\pi
t}}$, d'où : $$ L(\beta, x,t) = |\gamma| x K(\beta, \gamma,0,t)$$
Si $\beta = \gamma = 0$, $K(\beta,\gamma,0,t) = 1 \geq
\frac{1}{\sqrt{t}} \geq \sqrt{\frac{2}{\pi t}}$, et : $$ L(\beta, x,t)
\leq x K(\beta, \gamma, 0,t)$$ 
Si $\beta = 0$ et $\gamma > 0$, $K(\beta,\gamma,0,t)  \geq
\frac{1}{\gamma-\beta} \sqrt{\frac{2}{\pi t}}$, d'où : $$ L(\beta,x,t)
\leq x (\gamma - \beta) K(\beta,\gamma,0,t)$$ 
Si $\beta  > 0$ et $\gamma < \beta$, $L(\beta, x,t) =
\sqrt{\frac{2}{\pi t^3}} \frac{x}{\beta^2} + 2 \sinh (\beta x) e^{t
  \beta^2/2}$ et $K(\beta,\gamma,0,t) \geq \frac{1}{\beta-\gamma}
\sqrt{\frac{2}{\pi t^3}} + \frac{2 \beta}{\beta-\gamma} e^{t
  \beta^2/2}$. \\ \\
On en déduit que : $$ L(\beta,x,t) \leq \max
\left(\frac{x(\beta-\gamma)}{\beta^2},\sinh(\beta x)
  \frac{\beta-\gamma}{\beta} \right) K(\beta,\gamma,0,t)$$  
Or $x \leq \frac{\sinh(\beta x)}{\beta}$, d'où : $$ L(\beta,x,t) \leq
\max \left( \frac{\beta - \gamma}{\beta^3}, \frac{\beta -
    \gamma}{\beta} \right) \sinh(\beta x) K(\beta, \gamma, 0,t)$$
Si $\beta > 0$ et $\gamma = \beta$, on a $K(\beta,\gamma,0,t) \geq
\beta \sqrt{\frac{2}{\pi t^3}} + 2 e^{t \beta^2/2}$. \\ \\ On obtient
donc : $$L(\beta,x,t) \leq \max \left(\frac{x}{\beta^3}, \sinh(\beta
  x) \right) K(\beta, \gamma,0,t)$$ $$ \leq \max \left(\frac{1}{\beta^4} ,1
\right) \sinh(\beta x) K(\beta,\gamma,0,t)$$
Si $\beta > 0$ et $\gamma > \beta$, on a $K(\beta,\gamma,0,t) \geq
\frac{1}{\gamma-\beta} \sqrt{\frac{2}{\pi t^3}} + \frac{2
  \gamma}{\gamma - \beta}e^{t \beta^2/2}$ d'où : $$ L(\beta, \gamma,x)
\leq \max \left(x \frac{\gamma-\beta}{\beta^2}, \sinh(\beta x)
  \frac{\gamma-\beta}{\gamma} \right) K(\beta,\gamma,0,t) $$ $$\leq \max
\left(\frac{\gamma - \beta}{\beta^3}, \frac{\gamma - \beta}{\gamma}
\right) \sinh(\beta x) K(\beta,\gamma,0,t)$$
 \\ \\
Tout cela prouve que pour tous $\beta$, $\gamma$, il existe
$D(\beta,\gamma)$ tel qu'on ait, pour tout $t \geq 1$ et tout $x \geq
0$ : $$ L(\beta,x,t) \leq D(\beta,\gamma) \sinh((\beta^{+} + 1) x)
K(\beta,\gamma,0,t)$$
Maintenant, posons $\delta(\alpha) = \max \{ \alpha_m^{+}, m \in E
\}$, $\nu_{(E,\mu)} = \min \{ \mu_m, m \in E \}$ ($\nu_{(E,\mu)} > 0$ puisque $\mu_m
> 0$ pour tout $m \in E$), et $D(\alpha, \gamma) = \max \{ D(\alpha_m,
\gamma), m \in E \}$. \\ \\
On obtient les inégalités : $$ L(\alpha_k,x,t) \leq D(\alpha_k,
\gamma) \sinh((\alpha_k^+ + 1) x) K(\alpha_k,\gamma,0,t) $$ $$\leq
D(\alpha,\gamma) \sinh ((\delta(\alpha) + 1) x) \underset{m \in
  E}{\sum} \frac{\mu_m}{\nu_{(E,\mu)}} K(\alpha_m, \gamma,0,t)$$ $$
\leq \frac{D(\alpha, \gamma)}{\nu_{(E,\mu)}} \sinh((\delta(\alpha) +
1) x) Q_{(E,\mu)} (\alpha,\gamma, 0, 0,t)$$ 
On en déduit : $$ Q_{(E,\mu)} (\alpha, \gamma,x,k,t) \leq
\left(\frac{D(\alpha, \gamma)}{\nu_{(E,\mu)}} + C_{(E,\mu)}(\alpha,\gamma) \right)
\exp[(\delta(\alpha) + 1)x] Q_{(E,\mu)}(\alpha, \gamma,0,0,t) $$ inégalité
valable dès que $t \geq 1$, et que nous noterons : $$ Q_{(E,\mu)} (\alpha,
\gamma,x,k,t) \leq H_{(E,\mu)}(\alpha, \gamma) \exp(\psi(\alpha) x) Q_{(E,\mu)}(\alpha,
\gamma,0,0,t)$$ 
Cette inégalité nous permet de démontrer l'existence de la mesure
cherchée. En effet, si $s \geq 0$ et $\Gamma_s \in \mathcal{F}_s$, on a,
pour tout $t \geq s+1$ : 
$$ \mathbf{W}_{(E,\mu)}^{(t,\alpha,\gamma)} [\Gamma_s] = \mathbf{W}_{(E, \mu, 0,0)}
\left[ \mathbf{1}_{\Gamma_s} \frac{e^{\alpha_{N_t} X_t + \gamma
 L_t}}{\mathbf{W}_{(E, \mu, 0,0)}[e^{\alpha_{N_t} X_t + \gamma L_t} ]}
\right]$$ $$ = \mathbf{W}_{(E,\mu,0,0)} \left[ \mathbf{1}_{\Gamma_s}
  \frac{\mathbf{W}_ {(E, \mu, 0,0)}[e^{\alpha_{N_t} X_t + \gamma L_t} |
    \mathcal{F}_s ]}{\mathbf{W}_{(E, \mu, 0,0)} [e^{\alpha_{N_t} X_t +
      \gamma L_t} ]} \right]$$ $$ = \mathbf{W}_ {(E, \mu, 0,0)}\left[
  \mathbf{1}_{\Gamma_s} e^{\gamma L_s} \frac{\mathbf{W}_{(E, \mu, 0,0)}
    [e^{\alpha_{N_t} X_t + \gamma(L_t - L_s)}|X_s,N_s
    ]}{\mathbf{W}_{(E, \mu, 0,0)} [e^{\alpha_{N_t} X_t + \gamma L_t}]}
\right]$$ $$ = \mathbf{W}_ {(E, \mu, 0,0)}\left[\mathbf{1}_{\Gamma_s}
  \exp(\gamma L_s) \frac{R_{(E,\mu)}(\alpha, \gamma,X_s,N_s,t-s)}{R_{(E,\mu)}(\alpha,
    \gamma,0,0,t)} \right]$$ 
où $R_{(E,\mu)}(\alpha, \gamma, x,k,u) = \mathbf{W}_{(E, \mu, x,k)}
[\exp(\alpha_{N_u} X_u + \gamma L_u)]$. \\ \\
On sait que $\exp(\gamma L_s) \frac{R_{(E,\mu)} (\alpha, \gamma, X_s,N_s,
  t-s)}{R_{(E,\mu)}(\alpha, \gamma, 0,0,t)}$ est équivalent à $\exp(\gamma L_s)
\frac{Q_{(E,\mu)}(\alpha, \gamma,X_s,N_s,t-s)}{Q_{(E,\mu)} (\alpha, \gamma,0,0,t)}$
quand $t$ tend vers l'infini ($L_s$, $X_s$, $N_s$ étant fixés). \\ \\
Or $Q_{(E,\mu)}(\alpha, \gamma,x,k,u) = \underset{m \in E}{\sum} \mu_m
K(\alpha_m, \gamma,x,u) + L(\alpha_k,x,u)$ pour tous $x$, $k$, $u$,
donc d'après les estimations précédentes de $K$ et $L$, on a les
équivalents suivants : 
\\ \\ \\
\begin{tabular}{|p{3.5 cm}|p{7.5 cm}|} 
\hline 
Conditions sur $\alpha$, $\gamma$ & Equivalent de $Q_{(E,\mu)}(\alpha, \gamma,
x, k, u)$ pour $u \rightarrow \infty$ \\
\hline 
$\gamma \geq \alpha_m$ pour tout $m$, $\gamma = \alpha_m$ ssi $m \in
J$, $J$ sous-ensemble non vide de $E$, et $\gamma > 0$ & $ 2 \left(
  \underset{m \in J}{\sum} \mu_m \right) u \gamma^2 e^{-\gamma x + u
  \gamma^2/2}$ \\
\hline
$\gamma > \alpha_m$ pour tout $m \in E$ et $\gamma > 0$ & $\left( \underset{m
  \in E}{\sum} \frac{2 \gamma \mu_m}{\gamma- \alpha_m} \right) e^{-\gamma x
+ u \gamma^2/2}$ \\
\hline
$\alpha_m = \max(\alpha) = \bar{\alpha}$ ssi $m \in J$ ($J$
sous-ensemble non vide de $E$), $\bar{\alpha} > \gamma$ et
$\bar{\alpha} > 0$ & $e^{u \bar{\alpha}^2/2} \left( \frac{2
    \bar{\alpha}}{\bar{\alpha} - \gamma} \left( \underset{m \in
      J}{\sum} \mu_m \right) e^{- \bar{\alpha} x} + 2 \sinh
  (\bar{\alpha} x) \mathbf{1}_{k \in J} \right)$ \\
\hline
$\gamma = 0$, $\alpha_m = 0$ si $m \in J$ (sous-ensemble non vide de
$E$) et $\alpha_m < 0$ sinon & $\underset{m \in J}{\sum} \mu_m$ \\
\hline
$\gamma = 0$, $\alpha_m < 0$ pour tout $m \in E$ & $\sqrt{\frac{2}{\pi
    u}} \left( \underset{m \in E}{\sum} \frac{\mu_m}{|\alpha_m|}
\right)$ \\
\hline
$\alpha_m = 0$ si $m \in J$ (sous-ensemble non vide de $E$), $\alpha_m
< 0$ sinon, et $\gamma < 0$ & $\sqrt{\frac{2}{\pi u}} \left(
  \frac{\underset{m \in J}{\sum} \mu_m}{|\gamma|} + x \mathbf{1}_{k
    \in J} \right)$ \\
\hline
$\alpha_m < 0$ pour tout $m \in E$ et $\gamma < 0$ &
$\sqrt{\frac{2}{\pi u^3}} \left( \underset{m \in E}{\sum} \mu_m
  \frac{|\alpha_m| + |\gamma|}{\alpha_m^2 \gamma^2} + x \left(
    \frac{1}{\alpha_k^2} + \underset{m \in E}{\sum}
    \frac{\mu_m}{\alpha_m \gamma} \right) \right)$ \\
\hline
\end{tabular}
\\ \\ \\
On en déduit aisément que l'expression $\exp(\gamma L_s) \frac{R_{(E,\mu)}
  (\alpha, \gamma, X_s,N_s,t-s)}{R_{(E,\mu)}(\alpha, \gamma,0,0,t)}$ converge,
quand $t$ tend vers l'infini (à $L_s$, $X_s$, $N_s$ fixés) vers
$M_{(E,\mu)}(\alpha, \gamma, X_s, N_s, L_s)$, donné par le tableau suivant :  
\\ \\ \\
\begin{tabular}{|p{3.5 cm}|p{7.5 cm}|}
\hline
Conditions sur $\alpha$, $\gamma$ & $M_{(E,\mu)} (\alpha,\gamma,X_s,N_s,L_s)$
\\
\hline
$\gamma \geq \alpha_m$ pour tout $m$ et $\gamma > 0$ & $e^{\gamma(L_s
  - X_s) - s \gamma^2/2}$ \\
\hline 
$\alpha_m = \max(\alpha) = \bar{\alpha}$ ssi $m \in J$ ($J \subset E$
et $J \neq \emptyset$), $\bar{\alpha} > \gamma$ et $\bar{\alpha} > 0$ &
$e^{\gamma L_s - s \bar{\alpha}^2/2} \left( e^{- \bar{\alpha} X_s} +
  \frac{\bar{\alpha} - \gamma}{\bar{\alpha} \underset{m \in J}{\sum}
    \mu_m} \sinh(\bar{\alpha} X_s) \mathbf{1}_{N_s \in J} \right)$ \\
\hline
$\gamma = 0$, $\alpha_m \leq 0$ pour tout $m \in E$ & 1 \\
\hline
$\alpha_m = 0$ si $m \in J$ ($J \subset E$ et $J \neq \emptyset$),
$\alpha_m < 0$ sinon, et $\gamma < 0$ & $e^{\gamma L_s} \left( 1 +
  \frac{|\gamma|}{\underset{m \in J}{\sum} \mu_m} X_s \mathbf{1}_{N_s
    \in J} \right)$ \\
\hline
$\alpha_m < 0$ pour tout $m \in E$ et $\gamma <0$ & $e^{\gamma L_s}
\left( 1 + \frac{\frac{1}{\alpha_{N_s}^2} + \underset{m \in E}{\sum}
    \frac{\mu_m}{\alpha_m \gamma}}{\underset{m \in E}{\sum} \mu_m
    \frac{|\alpha_m| + |\gamma|}{\alpha_m^2 \gamma^2}} X_s \right)$
\\
\hline
\end{tabular}
\\ \\ \\
Par ailleurs, comme $t-s \geq 1$, on a les inégalités : $$R_{(E,\mu)}(\alpha,
\gamma, X_s, N_s, t-s) \leq Q_{(E,\mu)} (\alpha, \gamma, X_s, N_s, t-s) $$ $$\leq
H_{(E,\mu)}(\alpha, \gamma) e^{\psi(\alpha) X_s} Q_{(E,\mu)}(\alpha, \gamma,
0,0,t-s)$$
et $$R_{(E,\mu)} (\alpha, \gamma, 0,0,t) \geq \frac{1}{2} Q_{(E,\mu)}(\alpha,
\gamma,0,0,t)$$ pour $t$ assez grand (à $\alpha$, $\gamma$, $E$, $\mu$
fixés), puisque $R_{(E,\mu)}(\alpha, \gamma, 0,0,t)$ est équivalent à
$Q_{(E,\mu)}(\alpha, \gamma, 0,0,t)$ quand $t$ tend vers l'infini. \\ \\
De plus, pour $t$ assez grand : $$ \frac{Q_{(E,\mu)}(\alpha, \gamma,
  0,0,t-s)}{Q_{(E,\mu)}(\alpha, \gamma, 0,0,t)} \leq 2 M_{(E,\mu)} (\alpha,
\gamma,0,0,0) \leq 2$$
On en déduit que pour $t$ supérieur à une valeur ne dépendant que de
$E$, $\mu$, $\alpha$, $\gamma$ et $s$, on a : $$ e^{\gamma L_s}
\frac{R_{(E,\mu)}   (\alpha, \gamma, X_s, N_s, t-s)}{R_{(E,\mu)}(\alpha, \gamma, 0,0,t)} \leq 4
H_{(E,\mu)}(\alpha, \gamma) \exp(\psi(\alpha) X_s + \gamma L_s)$$
Ce majorant étant intégrable sous $\mathbf{W}_{(E, \mu, 0,0)}$, on en
déduit, par le théorème de convergence dominée : $$ \mathbf{W}_{(E,\mu)}^{(t, 
  \alpha, \gamma)} (\Gamma_s) \underset{t \rightarrow
  \infty}{\rightarrow} \mathbf{W}_{(E, \mu, 0,0)} [\mathbf{1}_{\Gamma_s}
M_{(E,\mu)} (\alpha, \gamma, X_s, N_s, L_s)]$$ 
On a donc prouvé l'existence d'une mesure $\mathbf{W}_{(E,\mu)}^{(\infty,
  \alpha, \gamma)}$ vérifiant : $$ \mathbf{W}_{(E,\mu)}^{(t, \alpha,
  \gamma)}(\Gamma_s) \underset{t \rightarrow \infty}{\rightarrow}
\mathbf{W}_{(E,\mu)}^{(\infty, \alpha, \gamma)}(\Gamma_s)$$ dès que $\Gamma_s \in
\mathcal{F}_s$ avec $s \geq 0$. \\ \\
De plus, il existe, sous $\mathbf{W}_{(E, \mu, 0,0)}$, une martingale
$(M_s^{(E,\mu,\alpha,\gamma)})_{s \geq 0}$, notée plus simplement
$(M_s)_{s \geq 0}$, égale à $M_{(E,\mu)}(\alpha, \gamma, X_s, N_s,L_s)$ et
vérifiant (pour $\Gamma_s \in \mathcal{F}_s$) : $$\mathbf{W}_{(E,\mu)}^{(\infty,
  \alpha, \gamma)} (\Gamma_s) = \mathbf{W}_{(E, \mu, 0,0)}
[\mathbf{1}_{\Gamma_s} M_s^{(E,\mu, \alpha, \gamma)}]$$
Ces dernières propriétés correspondent exactement à l'énoncé du théorème 1,
qui vient donc d'être démontré. 
\section{Etude du processus associé à $\mathbf{W}_{(E,\mu)}^{(\infty, \alpha,
    \gamma)}$} 
L'étude du processus associé à $\mathbf{W}_{(E,\mu)}^{(\infty, \alpha,
  \gamma)}$ se sépare en plusieurs cas, selon l'expression de la
martingale $(M_s)_{s \geq 0}$ précédemment donnée; ces cas
correspondent à la
 distinction effectuée dans l'énoncé du théorème 2. \\ \\ \\
\textbf{Cas où $\gamma \geq \alpha_m$ pour tout $m$ et $\gamma > 0$}
\\ \\
Sous $\mathbf{W}_{(E, \mu, 0,0)}$, $(Z_s = L_s - X_s)_{s \geq 0}$ est un
mouvement brownien. La densité de la loi de $(Z_u)_{0 \leq u \leq s}$
sous $\mathbf{W}_{(E,\mu)}^{(\infty, \alpha, \gamma)}$, par rapport à celle
d'un mouvement brownien sur $[0,s]$, est donc égale à $M_s = \exp(\gamma
  Z_s - s \gamma^2/2)$. \\ \\
$(Z_s)_{s \geq 0}$ est alors un mouvement brownien avec drift
$\gamma$ sous $\mathbf{W}_{(E,\mu)}^{(\infty, \alpha, \gamma)}$, et on retrouve $(X_s)_{s \geq 0}$ à partir de $(Z_s)_{s \geq
  0}$ grâce à l'expression : $X_s = \left( \underset{u \in
    [0,s]}{\sup} Z_u \right) - Z_s$. \\ \\
Autrement dit, $(X_s)_{s \geq 0}$ est la valeur absolue d'un processus bang-bang de
paramètre $\gamma$. \\ \\ 
Par ailleurs, $N_s$ n'intervient pas dans l'expression de $M_s$, donc
le processus $(N_s)_{s \geq 0}$, conditionné à $(X_s)_{s \geq 0}$, est
obtenu de la même manière sous $\mathbf{W}_{(E,\mu)}^{(\infty, \alpha,
  \gamma)}$ que sous $\mathbf{W}_{(E, \mu, 0,0)}$ : on choisit $(N_s)_{s \in
  I}$ pour chaque intervalle $I$ d'excursion de $(X_s)_{s \geq 0}$,
indépendamment et avec la loi $\mu$. \\ \\ \\
\textbf{Cas où $\alpha_m = \max(\alpha) = \bar{\alpha}$ ssi $m \in
  J$ ($J \subset E$ et $J \neq \emptyset$), $\bar{\alpha} > \gamma$
  et $\bar{\alpha} > 0$}  \\ \\
Avant de traiter ce deuxième cas en général, nous allons tout d'abord
supposer $J = \{ m \}$ pour un $m$ dans $E$, et $\gamma = 0$. \\
On a alors : $$M_s = e^{-s \bar{\alpha}^2/2} \left( e^{- \bar{\alpha}
    X_s} + \frac{1}{\mu_m} \sinh(\bar{\alpha} X_s) \mathbf{1}_{N_s =
    m} \right)$$ pour tout $s \geq 0$. \\ \\
Considérons à présent un processus $(Y_t, R_t)_{t \geq 0}$ sur
$\mathbf{R}_E$ défini de la manière suivante : \\ \\
- $(Y_t)_{t \geq 0}$ est la valeur absolue d'un mouvement brownien
avec drift $\bar{\alpha}$. \\ \\
- Soit $\mathcal{I}$ l'ensemble des intervalles d'excursion de
$(Y_t)_{t \geq 0}$. $R_t$ est alors défini comme étant constant sur
chaque intervalle $I \in \mathcal{I}$ ($R_t = R_I$ pour $t \in I$),
avec $R_I = m$ p.s. si $I$ est l'unique intervalle d'excursion non
borné de $(Y_t)_{t \geq 0}$, et avec les autres $(R_I)_{I \in
  \mathcal{I}}$ indépendants et de loi $\mu$. \\ \\
Montrons alors que $\mathbf{W}_{(E,\mu)}^{(\infty, \alpha, 0)}$ est la loi du
processus $(Y_t,R_t)_{t \geq 0}$. \\ \\
Pour cela, observons que la loi de $(Y_s,R_s)_{s \leq t}$,
conditionnellement à $(Y_t,R_t) = (x,k) \in E$, est égale à la loi de
$(X_s, N_s)_{s \leq t}$ sous $\mathbf{W}_{(E,\mu)}^{(\infty, \alpha, 0)}$,
conditionnellement à $(X_t, N_t) = (x,k)$. \\
Ces deux lois conditionnelles sont en effet égales à la loi du
processus $(Z_s,S_s)_{s \leq t}$ défini de la manière suivante : \\ \\
- $(Z_s)_{s \leq t}$ est la valeur absolue d'un pont brownien de $0$
vers $x$. \\ \\
- $(S_s)_{s \leq t}$ est constant sur les intervalles d'excursions de
$(Z_s)_{s \leq t}$; si on note $\mathcal{I}$ l'ensemble de ces
intervalles, et $S_s = S_I$ pour $s \in I$ et $I \in \mathcal{I}$,
alors les variables $(S_I)_{I \in \mathcal{I}}$ sont indépendantes et de loi
$\mu$, sauf pour l'intervalle $I$ de la forme $[a,t]$ ($a \in [0,t]$),
pour lequel on a $S_I = k$ p.s. \\ \\
L'égalité des lois conditionnelles précédentes (assez simple à
démontrer), entraîne que si pour tout $t$, la loi de $(Y_t,R_t)$ est
égale à celle de $(X_t, N_t)$ sous $\mathbf{W}_{(E,\mu)}^{(\infty, \alpha,
  0)}$, alors le processus $(Y_t, R_t)_{t \geq 0}$ a pour loi
$\mathbf{W}_{(E,\mu)}^{(\infty, \alpha, 0)}$, comme annoncé. \\ \\
Effectuons donc le calcul de la loi de $(Y_t,R_t)$ : $$ \mathbf{P}
(Y_t \in dy, R_t = k) = \mathbf{E} [\mathbf{1}_{Y_t \in dy} \mathbf{P}
(R_t = k | (Y_s)_{s \in \mathbf{R}_+})] $$ $$ = \mathbf{E}
[\mathbf{1}_{Y_t \in dy} (\mathbf{1}_{k=m} \mathbf{1}_{\forall s \geq
  t, Y_s > 0} + \mu_k \mathbf{1}_{\exists s \geq t, Y_s = 0})]$$
$$ = \mathbf{P} (Y_t \in dy) (\mathbf{1}_{k = m} \mathbf{P} (\forall s
\geq t, Y_s >0 | Y_t = y) + \mu_k \mathbf{P} (\exists s \geq t, Y_s =
0 | Y_t = y))$$
Rappelons que $Y_t= |B_t^{(\bar{\alpha})}|$ où
$(B_t^{(\bar{\alpha})})_{t \geq 0}$ est un mouvement brownien avec
drift $\bar{\alpha}$. On a donc : $$ \mathbf{P} (\exists s \geq t, Y_s
= 0 | Y_t=y) $$ $$= \mathbf{P} (\exists s \geq t, B_s^{(\bar{\alpha})} = 0
|B_t^{(\bar{\alpha})} = y) \mathbf{P} (B_t^{(\bar{\alpha})} = y |
|B_t^{(\bar{\alpha})}| = y) $$ $$ + \mathbf{P} (\exists s \geq t,
B_s^{(\bar{\alpha})} = 0 | B_t^{(\bar{\alpha})} = -y)
\mathbf{P}(B_t^{(\bar{\alpha})} = - y | |B_t^{(\bar{\alpha})}|= y) $$ $$
= e^{-2 \bar{\alpha} y}. \frac{e^{\bar{\alpha} y}}{e^{\bar{\alpha} y} +
  e^{- \bar{\alpha} y}} + 1. \frac{e^{- \bar{\alpha}
    y}}{e^{\bar{\alpha} y} + e^{- \bar{\alpha} y}} =
\frac{e^{-\bar{\alpha} y}}{\cosh(\bar{\alpha} y)}$$ et $$ \mathbf{P}
(\forall s \geq t, Y_s > 0 | Y_t = y) = 1 - \frac{e^{- \bar{\alpha}
    y}}{\cosh (\bar{\alpha} y)} = \frac{\sinh(\bar{\alpha}
  y)}{\cosh(\bar{\alpha} y)}$$
On en déduit : $$\mathbf{P}(Y_t \in dy, R_t = k) = \mathbf{P}(Y_t \in
dy) \left(\mathbf{1}_{k=m} \frac{\sinh(\bar{\alpha}
  y)}{\cosh(\bar{\alpha} y)} + \mu_k \frac{e^{- \bar{\alpha}
    y}}{\cosh(\bar{\alpha} y)} \right)$$ $$ = \mathbf{W}_{(E, \mu, 0,0)}
(X_t \in dy) \cosh (\bar{\alpha} y) e^{- t \bar{\alpha}^2/2} \left(
  \mathbf{1}_{k=m} \frac{\sinh(\bar{\alpha} y)}{\cosh(\bar{\alpha} y)}
  + \mu_k \frac{e^{-\bar{\alpha} y}}{\cosh(\bar{\alpha} y)} \right)$$
$$ = e^{-t \bar{\alpha}^2/2 } \mathbf{W}_{(E, \mu, x,k)} (X_t \in dy, N_t =k
) \left( e^{-\bar{\alpha} y} + \frac{1}{\mu_m} \sinh(\bar{\alpha} y)
  \mathbf{1}_{k=m} \right)$$
L'égalité des lois est donc démontrée. \\ \\
Ainsi, nous avons traité le cas particulier où $J = \{ m \}$ ($m \in
E$) et $\gamma = 0$. \\ \\
On remarque que si $E = \{ -1, 1 \}$, $\mu_1= \mu_{-1} = 1/2$ et $m=1$, le processus
$(X_t N_t)_{t \geq 0}$, qui est un mouvement brownien sous
$\mathbf{W}_{(E, \mu, 0,0)}$, est un mouvement brownien avec drift
$\bar{\alpha}$ sous $\mathbf{W}_{(E,\mu)}^{(\infty, \alpha,0)}$. \\
Cela se vérifie aussi bien avec la martingale $(M_s)_{s \geq 0}$
qu'avec la description du processus $(R_t,Y_t)_{t \geq 0}$ que nous
avons donnée ensuite. \\ \\ \\
A présent, traitons le cas, plus général, où $J = \{ m\}$ mais où
$\gamma$ n'est plus nécessairement nul. \\ \\
Dans ces conditions : $$ M_s = \exp(\gamma L_s - s \bar{\alpha}^2/2)
\left( e^{- \bar{\alpha} X_s} + \frac{\bar{\alpha} -
    \gamma}{\bar{\alpha} \mu_m} \sinh( \bar{\alpha} X_s)
  \mathbf{1}_{N_s = m} \right)$$
D'après le cas particulier précédent, la loi sous
$\mathbf{W}_{(E,\mu)}^{(\infty, \alpha, 0)}$ du temps local $L_{\infty}$ de
$(X_t,N_t)_{t \geq 0}$ est la même que celle du temps local d'un
mouvement brownien avec drift $\bar{\alpha}$ : $L_{\infty}$ est une
variable exponentielle de paramètre $\bar{\alpha}$. \\ \\
On en déduit l'existence de la mesure de probabilité $\nu$, donnée par
: $$\nu = \frac{\bar{\alpha} - \gamma}{\bar{\alpha}} \exp(\gamma
  L_{\infty}). \mathbf{W}_{(E,\mu)}^{(\infty, \alpha,0)}$$ (et sous laquelle
$L_{\infty}$ est une variable exponentielle de paramètre $\bar{\alpha}
- \gamma$). \\ \\
Montrons que $\nu$ est exactement la mesure $\mathbf{W}_{(E,\mu)}^{(\infty,
  \alpha,\gamma)}$ que nous étudions, celle-ci étant donc absolument
continue par rapport à $\mathbf{W}_{(E,\mu)}^{(\infty, \alpha, 0)}$. \\ \\
Pour prouver ce résultat, fixons $s \geq 0$ et $\Gamma_s \in
\mathcal{F}_s$. On a : 
$$\nu(\Gamma_s) = \mathbf{W}_{(E,\mu)}^{(\infty, \alpha,0)} \left[
  \frac{\bar{\alpha} - \gamma}{\bar{\alpha}} \mathbf{1}_{\Gamma_s}
  e^{\gamma L_s} e^{\gamma (L_{\infty} - L_s)} \right] $$ $$=
  \mathbf{W}_{(E,\mu)}^{(\infty, \alpha, 0)} \left[ \frac{\bar{\alpha} -
  \gamma}{\bar{\alpha}} \mathbf{1}_{\Gamma_s} e^{\gamma L_s}
  \mathbf{W}_{(E,\mu)}^{(\infty, \alpha,0)} [e^{\gamma (L_{\infty} -
  L_s)}|\mathcal{F}_s] \right] $$ $$= \mathbf{W}_{(E, \mu, 0,0)}\left[
  \frac{\bar{\alpha} - \gamma}{\bar{\alpha}} \mathbf{1}_{\Gamma_s}
  e^{\gamma L_s - s \bar{\alpha}^2/2} \left( e^{- \bar{\alpha} X_s} +
  \frac{1}{\mu_m} \sinh(\bar{\alpha} X_s) \mathbf{1}_{N_s =m} \right)
  \mathbf{W}_{(E,\mu)}^{(\infty, \alpha, 0)} [e^{\gamma (L_{\infty} - L_s)} |
  \mathcal{F}_s ] \right]$$
On observe alors les faits suivants (valables sous
  $\mathbf{W}_{(E,\mu)}^{(\infty, \alpha,0)}$) : \\ \\
- Conditionnellement à $\mathcal{F}_s$ et au fait que $X_t$ ne
  s'annule pour aucun $t \geq s$, $L_{\infty} - L_s$ est nul. \\ \\
- Conditionnellement à $\mathcal{F}_s$ et au fait que $X_t$ s'annule
  pour au moins un $t \geq s$, $L_{\infty} - L_s$ est une variable
  exponentielle de paramètre $\bar{\alpha}$. \\ \\ \\
Calculons maintenant la probabilité de chacun de ces deux cas,
  conditionnellement à $\mathcal{F}_s$. \\ \\
Pour cela, posons $T = \inf \{t \geq s, X_t = 0 \}$ et $A \in
  \mathcal{F}_s$ ($T$ est un temps d'arrêt). Pour $t \geq s$, on a :
  $$ \mathbf{W}_{(E,\mu)}^{(\infty, \alpha,0)} [T \leq t, A] =
  \mathbf{W}_{(E, \mu, 0,0)} [M_t^{(E,\mu, \alpha,0)}. \mathbf{1}_{T \leq t}
  . \mathbf{1}_A]$$ puisque $\{ T \leq t \}$ et $A$ sont
  $\mathcal{F}_t$-mesurables. \\ \\
De plus, $ \{ T \leq t \}$ et $A$ sont également $\mathcal{F}_{T
  \wedge t}$-mesurables, donc d'après le théorème d'arrêt :
  $$\mathbf{W}_{(E,\mu)}^{(\infty, \alpha, 0)} [T \leq t, A] =
  \mathbf{W}_ {(E, \mu, 0,0)} [M_{T \wedge t}^{(E,\mu, \alpha, 0)} \mathbf{1}_{T
  \leq t} \mathbf{1}_A]$$ $$ = \mathbf{W}_{(E, \mu, 0,0)}[M_T^{(E,\mu, \alpha, 0)}
  \mathbf{1}_{T \leq t} \mathbf{1}_{A} ]$$
Par convergence monotone : $$ \mathbf{W}_{(E,\mu)}^{(\infty, \alpha,0)} [T <
  \infty,A] = \mathbf{W}_{(E, \mu, 0,0)} [M_T^{(E, \mu,\alpha, 0)} \mathbf{1}_{T
  < \infty} \mathbf{1}_{A} ]$$ $$ = \mathbf{W}_{(E, \mu, 0,0)} [M_T^{(E, \mu,
  \alpha, 0)} \mathbf{1}_A]$$ puisque $T < \infty$ p.s. sous
  $\mathbf{W}_{(E, \mu, 0,0)}$. \\ \\
On a donc : $$ \mathbf{W}_{(E,\mu)}^{(\infty, \alpha,0)} [T< \infty, A] =
  \mathbf{W}_{(E, \mu, 0,0)} [\mathbf{1}_A \mathbf{W}_{(E, \mu,
  0,0)}[M_T^{(E,\mu, 
  \alpha,0)} | \mathcal{F}_s]]$$ $$ = \mathbf{W}_{(E,\mu)}^{(\infty, \alpha,
  0)} \left[ \mathbf{1}_{A} \frac{\mathbf{W}_{(E, \mu, 0,0)}
  [M_T^{(E,\mu, \alpha,0)} | \mathcal{F}_s]}{M_s^{(E,\mu, \alpha, 0)}}
  \right]$$ $$ = \mathbf{W}_{(E,\mu)}^{(\infty, \alpha, 0)} \left[
  \mathbf{1}_{A} \frac{\mathbf{W}_{(E, \mu, 0,0)} [e^{-(T-s)
  \bar{\alpha}^2/2} | \mathcal{F}_s]}{e^{-\bar{\alpha} X_s} +
  \frac{1}{\mu_m} \sinh( \bar{\alpha} X_s) \mathbf{1}_{N_s = m}}
  \right]$$
Or, conditionnellement à $\mathcal{F}_s$, $T-s$ est le temps
  d'atteinte de zéro d'un mouvement brownien issu de $X_s$ et
  indépendant de $X_s$. On en déduit : $$ \mathbf{W}_{(E, \mu, 0,0)}
  [e^{-(T-s) \bar{\alpha}^2/2} | \mathcal{F}_s] = e^{- \bar{\alpha}
  X_s}$$ et $$ \mathbf{W}_{(E,\mu)}^{(\infty, \alpha,0)} [T< \infty, A] =
  \mathbf{W}_{(E,\mu)}^{(\infty, \alpha, 0)} \left[ \mathbf{1}_A \frac{e^{-
  \bar{\alpha} X_s}}{e^{- \bar{\alpha} X_s} + \frac{1}{\mu_m}
  \sinh(\bar{\alpha} X_s) \mathbf{1}_{N_s = m}} \right]$$ autrement
  dit : $$ \mathbf{W}_{(E,\mu)}^{(\infty, \alpha, 0)} [T<\infty |
  \mathcal{F}_s] = \frac{e^{-\bar{\alpha} X_s}}{e^{-\bar{\alpha} X_s}
  + \frac{1}{\mu_m} \sinh(\bar{\alpha} X_s) \mathbf{1}_{N_s = m}}$$
Ceci permet d'écrire : $$ \mathbf{W}_{(E,\mu)}^{(\infty, \alpha,0)}
  [e^{\gamma(L_{\infty}-L_s)} | \mathcal{F}_s] =
  \frac{\frac{\bar{\alpha}}{\bar{\alpha} - \gamma} e^{-\bar{\alpha}
  X_s} + \frac{1}{\mu_m} \sinh(\bar{\alpha} X_s) \mathbf{1}_{N_s =
  m}}{e^{- \bar{\alpha} X_s} + \frac{1}{\mu_m} \sinh(\bar{\alpha}
  X_s) \mathbf{1}_{N_s = m}}$$ compte tenu des lois conditionnelles de
  $L_{\infty} - L_s$ précédemment données. \\ \\
Il en résulte : $$ \nu(\Gamma_s) = \mathbf{W}_{(E, \mu, 0,0)} \left[
  \mathbf{1}_{\Gamma_s} e^{\gamma L_s - s \bar{\alpha} ^2/2} \left( e^{-
  \bar{\alpha} X_s} + \frac{\bar{\alpha} - \gamma}{\bar{\alpha}
  \mu_m} \sinh(\bar{\alpha} X_s) \mathbf{1}_{N_s = m} \right) \right]
  $$ $$= \mathbf{W}_ {(E, \mu, 0,0)}[\mathbf{1}_{\Gamma_s} M_s^{(E,\mu, \alpha,
  \gamma)}]$$
On a donc l'égalité cherchée : $$\nu = \mathbf{W}_{(E,\mu)}^{(\infty, \alpha,
  \gamma)}$$
Nous venons donc de traiter le cas où $J$ est un singleton. \\ \\
Le cas général est alors facile à étudier; en effet la loi de
  $(X_t,N_t)_{t \geq 0}$, dans le cas général, est une moyenne des
  lois précédemment données, avec une pondération
  $\frac{\mu_m}{\underset{k \in J}{\sum} \mu_k}$ pour chaque $m \in
  J$. \\ \\
Autrement dit, le processus canonique sous $\mathbf{W}_{(E,\mu)}^{(\infty,
  \alpha, \gamma)}$ se décrit de la même manière qu'avant, sauf que sa
  dernière excursion se situe sur une branche quelconque appartenant à
  $J$, choisie alétoirement à l'aide de la mesure $\mu$. \\ \\ \\
 \textbf{Cas où $\gamma < 0$ et $\alpha_m \leq 0$ pour tout $m \in E$}
  \\ \\
Dans ce cas, on a : $$ M_s = e^{\gamma L_s}(1 + \theta_{N_s} X_s)$$ où
  les $(\theta_k)_{k \in E}$, positifs, dépendant de $\alpha$, sont
  tels que : $$ \underset{k \in E}{\sum} \mu_k \theta_k = |\gamma|$$
Nous allons tout d'abord supposer que $\theta_k =
  \frac{|\gamma|}{\mu_m}$ si $k=m$, $m$ étant un élément de $E$, et
  $\theta_k = 0$ si $k \neq m$. \\ \\
On a, dans ces conditions : $$M_s = e^{\gamma L_s} \left( 1 +
  \frac{|\gamma|}{\mu_m} X_s \mathbf{1}_{N_s=m} \right)$$ 
Considérons alors des réels positifs $l$ et $s$, une variable
  aléatoire $\mathcal{F}_{\tau_l}$-mesurable bornée $Y$ ($\tau_l$
  étant l'inverse du temps local de $(X_t)_{t \geq 0}$ pris en $l$),
  et une fonction $F$ mesurable bornée de $\mathcal{C}([0,s],
  \mathbf{R}_E)$ vers $\mathbf{R}$. \\ \\
On a alors, lorsque $t \geq 0$ (en utilisant le théorème d'arrêt pour
  la deuxième égalité) : $$ \mathbf{W}_{(E,\mu)}^{(\infty, \alpha, \gamma)}
  \left[ \mathbf{1}_{\tau_l \leq t} Y F((X_{\tau_l + u}, N_{\tau_l +
  u})_{0 \leq u \leq s}) \right]$$ $$= \mathbf{W}_{(E, \mu, 0,0)} \left[
  M_{t+s}^{(E, \mu,\alpha, \gamma)} \mathbf{1}_{\tau_l \leq t} Y
  F((X_{\tau_l + u}, N_{\tau_l+ u})_{0 \leq u \leq s}) \right] $$ $$=
  \mathbf{W}_{(E, \mu, 0,0)} \left[ M_{(t+s) \wedge (\tau_l +
  s)}^{(E,\mu, 
  \alpha, \gamma)} \mathbf{1}_{\tau_l \leq t} Y F((X_{\tau_l + u},
  N_{\tau_l + u})_{0 \leq u \leq s}) \right]$$ $$=
  \mathbf{W}_{(E, \mu, 0,0)} \left[ M_{\tau_l + s}^{(E,\mu, \alpha, \gamma)}
  F((X_{\tau_l + u}, N_{\tau_l + u})_{0 \leq u \leq s} )
  \mathbf{1}_{\tau_l \leq t} Y \right]$$
La convergence monotone entraîne alors (compte tenu du fait que
  $\tau_l < \infty$ p.s. sous $\mathbf{W}_{(E, \mu, 0,0)}$) : $$
  \mathbf{W}_{(E,\mu)}^{(\infty, \alpha, \gamma)} \left[ \mathbf{1}_{\tau_l <
  \infty} Y F((X_{\tau_l + u}, N_{\tau_l + u})_{0 \leq u \leq s})
  \right]$$ $$e^{\gamma l} \mathbf{W}_{(E, \mu, 0,0)} \left[e^{\gamma
  (L_{\tau_l + s} - L_{\tau_l})} \left( 1 + \frac{|\gamma|}{\mu_m}
  X_{\tau_l + s} \mathbf{1}_{N_{\tau_l + s} = m} \right) F((X_{\tau_l
  + u}, N_{\tau_l + u})_{0 \leq u \leq s}) Y \right]$$
D'après la propriété de Markov de l'araignée, $(X_{\tau_l + u},
  N_{\tau_l + u})_{0 \leq u \leq s}$ est indépendant de
  $\mathcal{F}_{\tau_l}$ sous $\mathbf{W}_{(E, \mu, 0,0)}$ et a la même loi que
  $(X_u, N_u)_{0 \leq u \leq s}$. \\ \\
On en déduit facilement : $$ \mathbf{W}_{(E,\mu)}^{(\infty, \alpha, \gamma)}
  \left[ \mathbf{1}_{L_{\infty} \geq l} Y F((X_{\tau_l +u}, N_{\tau_l
  + u})_{0 \leq u \leq s}) \right]$$ $$ =e^{\gamma l}
  \mathbf{W}_{(E,\mu)}^{(\infty, \alpha, \gamma)} [F(X_u,N_u)_{0 \leq u \leq
  s}] \mathbf{W}_{(E, \mu, 0,0)} [Y]$$
En particulier, pour $F$ et $Y$ égaux à $1$, on obtient :
  $$\mathbf{W}_{(E,\mu)}^{(\infty, \alpha, \gamma)} [L_{\infty} \geq
  l] = \exp(\gamma l)$$
On a donc les caractéristiques suivantes : \\ \\
- $L_{\infty}$ est une variable exponentielle de paramètre
  $|\gamma|$. \\ \\
- Conditionnellement à $L_{\infty} \geq l$, $(X_s,N_s)_{0 \leq s \leq
  \tau_l}$ est une araignée brownienne arrêtée en $\tau_l$, et
  $(X_{\tau_l + s}, N_{\tau_l + s})_{s \geq 0}$ admet pour loi
  $\mathbf{W}_{(E,\mu)}^{(\infty, \alpha, \gamma)}$; de plus, ces deux
  processus sont indépendants. \\ \\
On déduit de ce qui précède que conditionnellement à $L_{\infty} = l$,
  $(X_s,N_s)_{0 \leq s \leq \tau_l}$ est encore une araignée arrêtée
  en $\tau_l$, et $(X_{\tau_l + s},N_{\tau_l + s})_{s \geq 0}$ est un
  processus de loi $\mathbf{W}_{(E,\mu)}^{(\infty, \alpha, \gamma)}$,
  conditionné par le fait qu'il ne s'annule qu'au temps zéro; les deux
  processus étant encore indépendants. \\ \\
Pour décrire le deuxième processus, considérons $s \geq 0$, $\Gamma_s \in
  \mathcal{F}_s$, $l \geq 0$ et $t \geq s$. On a : $$
  \mathbf{W}_{(E,\mu)}^{(\infty, \alpha, \gamma)} [\Gamma_s, \tau_l \leq t] =
  \mathbf{W}_{(E, \mu, 0,0)} [M_t^{(E,\mu, \alpha, \gamma)} \mathbf{1}_{\Gamma_s}
  \mathbf{1}_{\tau_l \leq t}]$$ $$ = \mathbf{W}_{(E, \mu, 0,0)} [M_{\tau_l \vee
  s}^{(E, \mu, \alpha, \gamma)} \mathbf{1}_{\Gamma_s} \mathbf{1}_{\tau_l \leq
  t}]$$ d'où $$\mathbf{W}_{(E,\mu)}^{(\infty, \alpha, \gamma)} [\Gamma_s, \tau_l
  < \infty] = \mathbf{W}_ {(E, \mu, 0,0)}[M_{\tau_l \vee s}^{(E,\mu, \alpha,
  \gamma)} \mathbf{1}_{\Gamma_s}]$$ et donc : $$\mathbf{W}_{(E,\mu)}^{(\infty,
  \alpha, \gamma)} [\Gamma_s, L_{\infty} \leq l] = \mathbf{W}_{(E, \mu, 0,0)}
  [(M_s^{(E,\mu, \alpha, \gamma)} - M_{\tau_l \vee s}^{(E,\mu, \alpha,
  \gamma)}) \mathbf{1}_{\Gamma_s} ] $$ $$ = \mathbf{W}_{(E, \mu, 0,0)} \left[
  \left( e^{\gamma L_s} \left( 1 + \frac{|\gamma|}{\mu_m} X_s
  \mathbf{1}_{N_s = m} \right) - e^{\gamma l} \right) \mathbf{1}_{L_s
  \leq l} \mathbf{1}_{\Gamma_s} \right] $$ $$ = \mathbf{W}_{(E, \mu, 0,0)}
  [L_s \leq l] \mathbf{W}_ {(E, \mu, 0,0)}\left[ \mathbf{1}_{\Gamma_s} \left(
  e^{\gamma L_s} \left( 1 + \frac{|\gamma|}{\mu_m} X_s
  \mathbf{1}_{N_s = m} \right) - e^{\gamma l} \right) | L_s \leq l
  \right]$$
Comme $\mathbf{W}_{(E,\mu)}^{(\infty, \alpha, \gamma)} [L_{\infty} \leq l] = 1
  - e^{\gamma l}$, on a : $$\mathbf{W}_{(E,\mu)}^{(\infty, \alpha, \gamma)}
  [\Gamma_s | L_{\infty} \leq l] $$ $$= \frac{\mathbf{W}_{(E, \mu, 0,0)} [L_s \leq
  l]}{1-e^{\gamma l}} \mathbf{W}_{(E, \mu, 0,0)}\left[ \mathbf{1}_{\Gamma_s}
  \left( e^{\gamma L_s} \left( 1 + \frac{|\gamma|}{\mu_m} X_s
  \mathbf{1}_{N_s = m} \right) - e^{\gamma l} \right) | L_s \leq l
  \right] $$ $$= \frac{\mathbf{W}_ {(E, \mu, 0,0)}[L_s \leq l]}{1 -
  e^{\gamma l}} \tilde{\mathbf{W}} (l) \left[ \mathbf{1}_{\Gamma_s}
  \left( e^{\gamma L_s} \left( 1 + \frac{|\gamma|}{\mu_m} X_s
  \mathbf{1}_{N_s = m} \right) - e^{\gamma l} \right) \right]$$ où
  $\tilde{\mathbf{W}} (l)$ désigne la loi de $(X_s,N_s)_{s \geq 0}$
  conditionnée par le fait que $L_s \leq l$. \\ \\
Quand $l$ tend vers zéro, $\frac{\mathbf{W}_{(E, \mu, 0,0)} [L_s \leq l]}{1
  - e^{\gamma l}}$ tend vers $\frac{1}{|\gamma|} \sqrt{\frac{2}{\pi
  s}}$. \\ \\
D'autre part, si $l$ et $L_s$ tendent vers zéro à $(X_s,N_s)$ fixé,
  $e^{\gamma L_s} \left( 1 + \frac{|\gamma|}{\mu_m} X_s
  \mathbf{1}_{N_s = m} \right) - e^{\gamma l}$ tend vers
  $\frac{|\gamma|}{\mu_m} X_s \mathbf{1}_{N_s = m}$. \\ \\
Ceci permet de démontrer : $$ \mathbf{W}_{(E,\mu)}^{(\infty, \alpha, \gamma)}
  [\Gamma_s | L_{\infty} = 0] = \tilde{\mathbf{W}} (0)
  \left[\mathbf{1}_{\Gamma_s} \sqrt{\frac{2}{\pi s}} \frac{X_s}{\mu_m}
  \mathbf{1}_{N_s = m} \right]$$ où $\tilde{\mathbf{W}} (0)$ est la
  loi d'une araignée sur $[0,s]$, conditionnée par sa non-annulation
  en dehors du temps $0$; on remarque que sous $\tilde{\mathbf{W}}
  (0)$, $(X_u)_{u \leq s}$ est un méandre brownien. \\ \\
On en déduit alors que sous $\mathbf{W}_{(E,\mu)}^{(\infty, \alpha, \gamma)}$,
  et conditionnellement au fait que $X_s > 0$ pour tout $s > 0$,
  $(X_s)_{s \geq 0}$ est un processus de Bessel de dimension 3, et
  $N_s =m$ pour tout $s$. \\ \\
On a donc la description de $(X_t,N_t)_{t \geq 0}$ dans le cas où un
  seul des $\theta_k$ précédemment donnés est nul. \\ \\
Le cas général est simple à étudier à présent; en effet, il suffit de
  faire une moyenne pondérée des mesures précédemment décrites pour
  chacun des $m \in E$ (avec la pondération $\frac{\mu_m
  \theta_m}{|\gamma|}$). \\ \\ \\
\textbf{Cas où $\gamma = 0$ et $\alpha_m \leq 0$ pour tout $m \in E$}
  \\ \\
Ce cas est le plus simple de tous : $\mathbf{W}_{(E,\mu)}^{(\infty, \alpha, 0)}$
  est exactement la loi d'une araignée brownienne, puisque $M_s^{(E,\mu,
  \alpha, 0)}$ est constante et égale à 1. \\ \\ \\
Nous avons maintenant étudié tous les cas possibles pour $\gamma$ et
  $\alpha$, et il est facile de vérifier que cette étude entraîne le
  théorème 2. \\ \\
De plus, si $\gamma < 0$ et $\alpha_m \leq 0$ pour tout $m$, on a
  (pour tout $m \in E$ et avec les notations du théorème 2) :
  $$\mathbf{P} (M = m) = \frac{\mu_m}{\underset{k \in J}{\sum}
  \mu_k} \mathbf{1}_{m
  \in J}$$ si $J= \{ m \in E, \alpha_m = 0 \}$ est non vide, et
  $$\mathbf{P} (M = m) =
  \frac{ \mu_m \left( \frac{|\gamma|}{\alpha_m^2} + \underset{k \in
  E}{\sum} \frac{\mu_k}{|\alpha_k|} \right)} {\underset{k \in
  E}{\sum} \mu_k \frac{|\alpha_k| + |\gamma|}{\alpha_k^2}}$$ si $J =
  \emptyset$. \\ \\
Par ailleurs, on peut observer que $(X_s,N_s)_{s \geq 0}$ s'annule pour des valeurs
  arbitrairement grandes de $s$ ssi $\gamma \geq 0$ et $\gamma \geq
  \alpha_m$ pour tout $m$. \\ \\ \\
\textbf{Remarque 1 : } Le théorème 2, que nous venons de prouver,
  indique différents comportements possibles pour le processus limite
  obtenu, selon les valeurs des réels $\alpha_m$ ($m \in E$) et
  $\gamma$. \\
Cette distinction de cas généralise celle que l'on obtient à partir
  des résultats démontrés dans [8]. \\ \\
Par ailleurs, on observe que la distinction de cas donnée dans [8] est
  étroitement liée à celle que Y. Harriya et M. Yor obtiennent dans
  [4]; ce lien peut vraisemblablement être expliqué en comparant le 
 comportement des quantités $\int_0^t \exp (2X_s) ds$ et
  $\exp(2S_t)$, $(X_s)_{s \geq 0}$ étant un mouvement brownien et
  $S_t$ son maximum sur $[0,t]$. \\ \\
On peut alors se demander s'il est possible de mettre en évidence des liens
  analogues entre des
  pénalisations d'araignées browniennes. \\ \\ 
\textbf{Remarque 2 : } Soit $\nu$ une mesure de probabilité définie sur
  $\mathcal{C} (\mathbf{R}_+, \mathbf{R}_E)$, dont la densité par
  rapport à $\mathbf{W}_{(E, \mu, 0,0)}$, conditionnellement à
  $\mathcal{F}_s$ ($s \geq 0$), existe et s'écrit sous la forme : $$
  g(s,X_s,N_s) = \exp(-s \beta^2/2) f_{N_s} (X_s)$$ avec $f_m \in
  \mathcal{C}^2 (\mathbf{R}_+)$ pour tout $m \in E$, $f_m(0)$ ne
  dépendant pas de $m$ (ce qui permet de poser $f_0(0) =f_m(0)$), et
  $\beta > 0$. \\ \\
$g(s,X_s,N_s)$ est une martingale sous $\mathbf{W}_{(E, \mu, 0,0)}$. \\
L'étude du générateur infinitésimal de l'araignée permet alors de montrer
  les faits suivants : \\ \\
- Pour tous $x \in \mathbf{R}_+^*$, $s \geq 0$ et $m \in E$,
  $\frac{\partial g}{\partial s} (s,x,m) + \frac{1}{2}
  \frac{\partial^2 g}{\partial x^2} (s,x,m) = 0$. \\ \\
- Pour tout $s \geq 0$, $\underset{m \in E}{\sum} \mu_m
  \frac{\partial g}{\partial x} (s,0,m) = 0$. \\ \\
La première égalité donne : $$ - \frac{\beta^2}{2} e^{-s \beta^2/2}
  f_m(x) + \frac{1}{2} e^{-s \beta^2/2} f''_m(x) = 0$$
soit $f''_m(x) = \beta^2 f_m (x)$. \\ \\
On en déduit qu'il existe $\delta_m$ et $\lambda_m \in \mathbf{R}$
  tels que : $$f_m(x) = \delta_m \exp(-\beta x) + \lambda_m \sinh (\beta
  x)$$ pour tout $x \geq 0$. \\ \\
De plus, comme $g(s,X_s,N_s)$ est une densité, son espérance sous
  $\mathbf{W}_{(E, \mu, 0,0)}$ est 1. En particulier, pour $s=0$, on obtient
  $g(0,0,0) = 1$, donc $f_0(0) = 1$, ce qui implique $f_m(0) = 1$ et
  $\delta_m = 1$. \\ \\
Par ailleurs, $f_m(x) \geq 0$ pour tout $x \geq 0$, donc $\lambda_m
  \geq 0$. \\ \\
La deuxième égalité à vérifier implique alors : $\underset{m \in
  E}{\sum} \mu_m f'_m(0) = 0$, soit $\underset{m \in E}{\sum} \mu_m
 (1 - \lambda_m) = 0$ et $\underset{m \in E}{\sum} \mu_m \lambda_m =
  1$. \\ \\
On en déduit que $\nu$ est une moyenne pondérée des mesures
  $\mathbf{W}_{(E,\mu)}^{(\infty, \alpha, 0)}$ obtenues en prenant
  successivement, pour chaque $m$, $\alpha_m = \beta > \alpha_k$ (pour
  tout $k \neq m$), la pondération étant $\mu_m \lambda_m$. \\ \\
Les processus obtenus sont en fait des généralisations du mouvement
  brownien avec drift. 
\section*{Bibliographie}
\noindent
[1] M. Barlow, J. Pitman, M. Yor : On Walsh's Brownian motions,
Séminaire de Probabilités XXIII, 275-293 (1989) \\
\noindent
[2] M. Barlow, J. Pitman, M. Yor : Une extension multidimensionnelle
de la loi de l'arc sinus, Séminaire de Probabilités XXIII, 294-314
(1989) \\
\noindent
[3] A.-S. Cherny, A.-N. Shiryaev : Some distributional properties of a
Brownian motion with a drift and an extension of P. Lévy's theorem,
SIAM Theory of Probability and Its Applications \textbf{44}, 412-418
(1999) \\
\noindent
[4] Y. Hariya, M. Yor : Limiting distributions associated with moments
of exponential Brownian functionals, Studia
Sci. Math. Hungar. \textbf{41}, 193-242 (2004) \\
\noindent 
[5] J. Pitman : The distributions of local times of a Brownian bridge,
Séminaire de Probabilités XXXIII, 388-394 (1999) \\
\noindent
[6] B. Roynette, P. Vallois, M. Yor : Limiting laws associated with
Brownian motion perturbated by normalized exponential weights,
C.R.A.S. Paris, Sér. I Math \textbf{337}, 667-673 (2003) \\
\noindent
[7] B. Roynette, P. Vallois, M. Yor : Limiting laws associated with
Brownian motion perturbated by its maximum, minimum and local time,
II, to appear in Studia Sci. Math. Hungar. (2005) \\
\noindent
[8] B. Roynette, P. Vallois, M. Yor : Limiting laws for long Brownian
bridges perturbed by their one-sided maximum, III, to appear in 
Periodica Hungarica (2005) \\
\noindent
[9] J.-B. Walsh : A diffusion with a discontinuous local time, Temps
Locaux, Astérisque \textbf{52-53}, 37-45 (1978)
\end{document}